\documentclass{article}

\usepackage{arxiv}

\usepackage[utf8]{inputenc} 
\usepackage[T1]{fontenc}    
\usepackage{hyperref}       
\usepackage{url}            
\usepackage{booktabs}       
\usepackage{amsfonts}       
\usepackage{nicefrac}       
\usepackage{microtype}      
\usepackage{lipsum}

\usepackage[final]{graphicx}
\usepackage{amsmath}
\usepackage{amsthm}

\newtheorem{theorem}{Theorem}

\newtheorem{definition}{Definition}

\hypersetup{unicode=true}
\usepackage[all]{xy}

\title{LEHMER PAIRS AND BINOMIAL SERIES
}

\author{
 Kirill V.~Kapitonets \\
 BAUMAN MSTU, GRADUATE 1990\\
  MCC EuroChem\\
  Moscow\\
  Russian Federation\\
  \texttt{kkapitonets@live.com} \\
}

\begin{document}
\maketitle
\begin{abstract}
The Hardy function $Z(t)=\zeta(1/2+it)e^{i\theta(t)}$ takes real values for real $t$ and its real zeros are zeros $\zeta(s)$ on the critical line $1/2+it$.
\par 
After discovering the critical value of the local maximum in 1956, Lehmer formulated the assumption that the Hardy function could have a negative local maximum or a positive local minimum.
\par 
In the paper the Generalized Hardy function is defined as the real part of the Hardy function on any line $\alpha_\nu+it$ parallel to the critical line $1/2+it$
$$Z_{\alpha_\nu}(t)=Re\  \zeta(\alpha_\nu+it)e^{i\theta(t)}$$
\par 
and established an distinct relationship between the zeros of the $\cos\theta(t)$ function and the zeros of the Generalized Hardy function.
\begin{equation}
\notag
\begin{gathered}
\forall \Delta T_\lambda=(t_\lambda, t_{\lambda+1}],\ 
t_\lambda=2\pi\lambda^2,\ \lambda=1,\ 2,\ 3\ ...\\
\exists A_\lambda:\forall \hat\alpha_\lambda>A_\lambda\\
|\cos\theta(t) -Z_{\hat\alpha_\lambda}(t)|<\epsilon(A_\lambda),\ t\in \Delta T_\lambda
\end{gathered}
\end{equation}
\par 
Then the binomial series is used to establish a relationship between the values of the Generalized Hardy function on any two lines $\alpha_\nu+it$ and $\alpha_{\nu+1}+it$ parallel to the critical line.
\begin{equation}
\notag
\begin{gathered}
Z_{\alpha^{(\lambda)}_\nu, 0}=Z_{\alpha^{(\lambda)}_{\nu+1}}\\
Z_{\alpha^{(\lambda)}_\nu, k}=Z_{\alpha^{(\lambda)}_\nu, k-1}-P_{\alpha^{(\lambda)}_\nu,k}
\end{gathered}
\end{equation}
\begin{equation}
\notag
\begin{split}
P_{\alpha^{(\lambda)}_\nu,k}=\frac{(-1)^k}{k!}\big\{\alpha^{(\lambda)}_{\nu+1}(\alpha^{(\lambda)}_{\nu+1}-1)(\alpha^{(\lambda)}_{\nu+1}-2)&...(\alpha^{(\lambda)}_{\nu+1}-k+1)\\
-\alpha^{(\lambda)}_{\nu}(\alpha^{(\lambda)}_{\nu}-1)(\alpha^{(\lambda)}_{\nu}-2)&...(\alpha^{(\lambda)}_{\nu}-k+1)\big\}Q_k
\end{split}
\end{equation}
\begin{equation}
\notag
\begin{gathered}
Q_k=\lim_{m\to\infty}\sum^{m+k_m}_{n=2}\Big(\frac{n-1}{n}\Big)^kT_n,\ k=0,1,2...\\
T_n=\cos(\theta(t)-t\log(n))\delta^{(\nu_m)}_{\lambda_{\omega},n}
\end{gathered}
\end{equation}

\par 
Thus, by induction between values $\sigma=1/2$ and $\sigma=\hat\alpha_\lambda>A_\lambda$
$$\alpha^{(\lambda)}_1<\alpha^{(\lambda)}_2<\alpha^{(\lambda)}_3<...<\alpha^{(\lambda)}_{\nu}<...<\alpha^{(\lambda)}_{\mu_\lambda}$$ 
an distinct relationship has been established between the zeros of the function $\cos\theta(t)$ and  the zeros of the Hardy function.
\end{abstract}

\section{Introduction}
In 1956, calculating the first 15,000 zeros of the Zeta function on a computer, Lehmer \cite{LD} discovered a pair of zeros whose ordinates differ by less than $0.1$, or rather, by $0.037698499$ and assumed that the Hardy \cite{AI1, VK} function can have a negative local maximum or a positive local minimum.
\begin{figure}[ht!]
\centering
\includegraphics[scale=0.5]{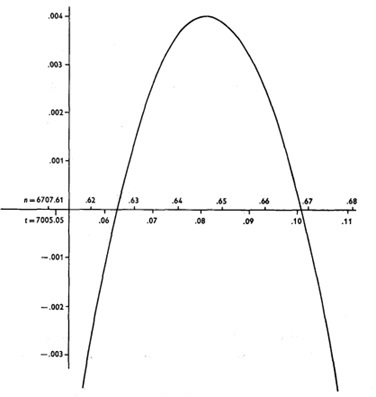}
\caption{Lehmer pair}
\label{fig:lehmer_pair}
\end{figure}
\par 
The Hardy function $Z(t)$ takes real values for real $t$ and its real zeros are zeros $\zeta(s)$ on the critical line $1/2+it$.
\begin{equation}
\label{hardy_function_1}
\begin{gathered}
Z(t)=\zeta(1/2+it)e^{i\theta(t)}\\
e^{i\theta(t)}=\pi^{-it/2}\frac{\Gamma(1/4+it/2)}{|\Gamma(1/4+it/2)|}\\
\zeta(s)=\sum^\infty_{n=1} n^{-s}\\
s=\sigma+it
\end{gathered}
\end{equation}
\par
\begin{figure}[ht!]
\centering
\includegraphics[scale=0.7]{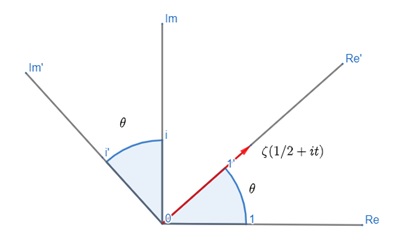}
\caption{Hardy function $\zeta(1/2+it)e^{i\theta(t)}$}
\label{fig:rotation_1}
\end{figure}

 
Obviously (Fig. \ref{fig:rotation_1}), to obtain the Hardy function, we rotate the coordinate axes of the complex plane by an angle of $\theta(t)$ to align the vector corresponding to $\zeta(1/2+it)$ with the real axis, i.e. $Z(t)$ becomes a real function.
\par 
The angle $\theta(t)$ occurs if the value $1/2+it$ is substituted into the functional equation \cite{TI} of the Zeta function
\begin{equation}\label{zeta_func_eq}\pi^{-s/2}\Gamma(\frac{s}{2})\zeta(s)=\pi^{-(1-s)/2}\Gamma(\frac{1-s}{2})\zeta(1-s);\end{equation}
\par 
In 2004, calculating the zeros of the Zeta function over a large interval, Gordon \cite{GX} provided a table of Lehmer pairs whose ordinates differ by less than $0.0001$.
\begin{center}
\includegraphics[scale=0.35]{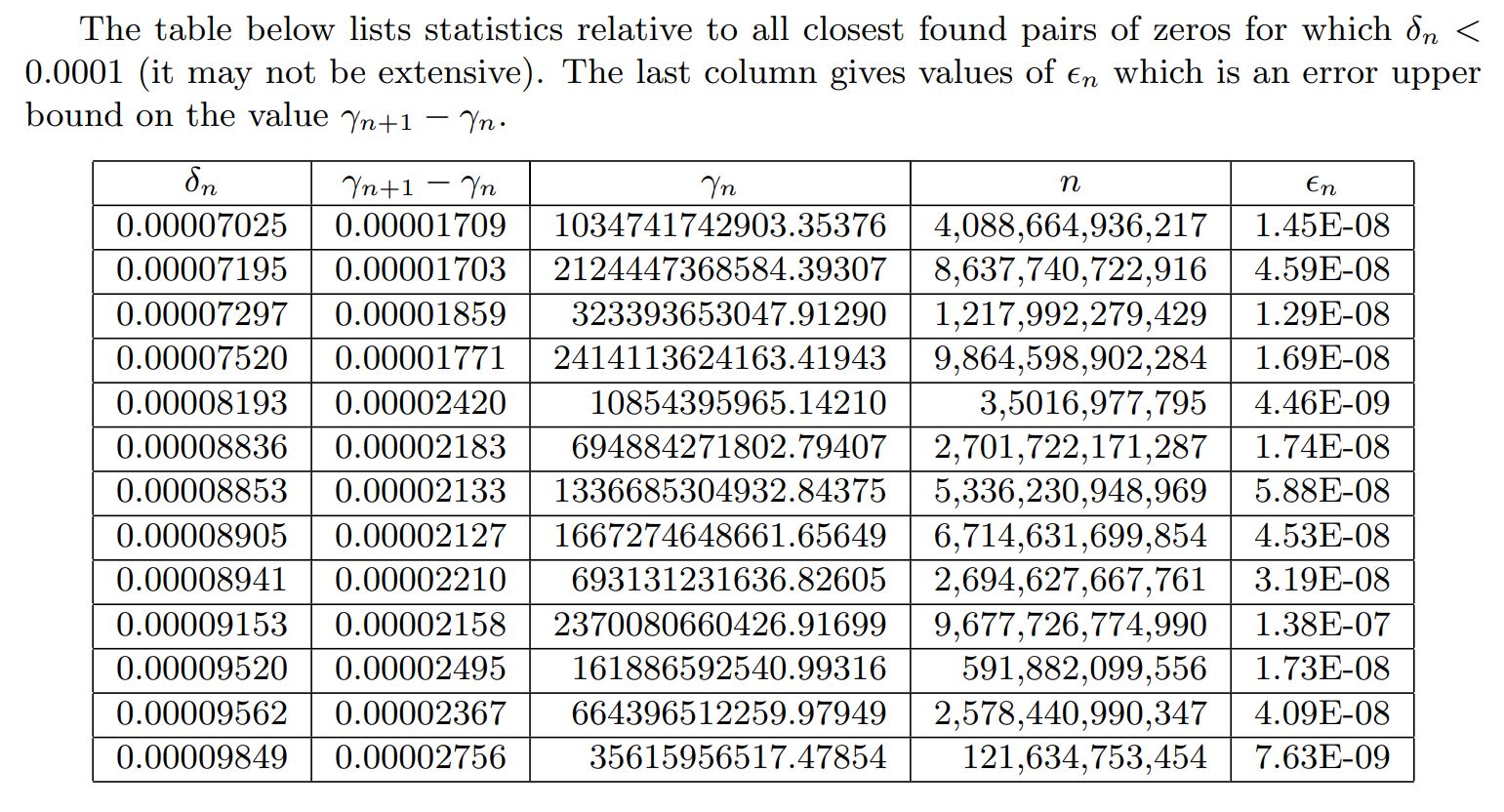}
\end{center}
\begin{figure}[ht!]
\centering
\caption{Lehmer pair list}
\label{fig:lehmer_pair_list}
\end{figure}
\par 
Such pairs of zeros are now called Lehmer pairs, and the Lehmer phenomenon has become one of the main arguments \cite{AI2} against the truth of the Riemann hypothesis.
\par 
In 1994, a group of researchers Chordas, Smith, and Varga \cite{CSV} proposed defining a Lehner pair through the zeros of the Zeta function to determine the lower bound of the de Bruijn-Newman constant (in Stopple \cite{ST} notation)
\begin{equation}
\notag
\begin{gathered}
\Delta^2\sum_{\gamma\ne\gamma_-, \gamma_+}\frac{1}{(\gamma-\gamma_-)^2}+\frac{1}{(\gamma-\gamma_+)^2}<\frac{4}{5}\\
\Delta=\gamma_+- \gamma_-
\end{gathered}
\end{equation}
\par 
This definition is based on the statistics of zeros of the Zeta function calculated at that time, i.e. up to ordinate 1,090,879,645.
\par 
In 2015, Stopple \cite{ST} proposed the definition of strict Lehmer pairs in terms of the pre-Schwarzian derivative of the Xi function
\begin{equation}
\notag
\begin{gathered}
\Delta^2\big(-P\Xi'(\gamma_+)-P\Xi'(\gamma_-)\big)\big)<\frac{42}{5}\\
\Delta=\gamma_+- \gamma_-\\
P\Xi'=\Big(\frac{\Xi''}{\Xi'}\Big)'\\
\Xi=\xi(1/2+it)\\
\xi(s)=\frac{s(1-s)}{2}\pi^{-s/2}\Gamma(\frac{s}{2})\zeta(s)
\end{gathered}
\end{equation}
which corresponds to the real part of the integral distance between the pair of zeros of the Xi function and the zeros of the derivative of the Xi function
\begin{equation}
\notag
-\Delta^2Re\Big\{\sum_{\substack{|\gamma'-\gamma_0'|\le1/\log\log(\gamma_0')\\ \rho'\ne\rho'_0}}
\Big(\frac{1}{(\rho_+-\rho')^2}+\frac{1}{(\rho_--\rho')^2}\Big)\Big\}
\end{equation}

\par 
Obviously, the definition of $\xi(s)$ is obtained from the functional equation (\ref{zeta_func_eq}) by adding the multiplier $\frac{s(1-s)}{2}$, and the Hardy function is related \cite{TI} with Xi function by the equation
\begin{equation}
\notag
Z(t)=-2\pi^{1/4}\frac{\Xi(t)}{(t^2+1/4)|\Gamma(1/4+it/2)|}
\end{equation}

\par 
In 2017, Simonic \cite{SI} proposed the definition of strict Lehmer pairs in terms of the pre-Schwarzian derivative of the Hardy function, just as Stopple does it
\begin{equation}
\notag
\begin{gathered}
\hat g_{\{\gamma_1,\gamma_2\}}=\frac{1}{3}(\gamma_1-\gamma_2)^2\Big(\hat F(\gamma_1)+\hat F(\gamma_2)\Big)-2<\frac{4}{5}-3(\gamma_1-\gamma_2)^2\Big(\frac{1}{\gamma_1^2}+\frac{1}{\gamma_2^2}\Big)\\
\hat F(t)=\frac{Z'''}{Z''}(t)+\frac{3}{4}\Big(\frac{Z''}{Z'}\Big)^2(t)
\end{gathered}
\end{equation}
\par 
In 2018, Rogers, Brad, and Tao \cite{RBT} prove a general relation for the zeros of the Xi function and do not use any definition of Lehmer pairs.
\begin{equation}
\notag
\log\frac{1}{|x_{j+1}-x_j|}\ll(\log^2j)\log\log j
\end{equation}

\par 
Thus, different researchers propose different integral estimates of the distance between the zeros of the Zeta function in order to define Lehmer pairs.
\par 
We review any neighboring pair of zeros $\{\gamma_n,\gamma_{n+1}\}$ of the Hardy function (\ref{hardy_function_1}) and prove the theorem
\begin{theorem}{on the transformation of zeros of the Hardy function}
\par 
i) the zeros of the Hardy function can be obtained by transforming the zeros of the function $\cos\theta(t)$
\begin{equation}
\notag
\begin{gathered}
\exists\{f_k\}^\infty_{k=1}:\forall \overline t_n,\ \overline t_{n+1}\in\{\overline t_n\}^\infty_{n=1}\\ 
\{\overline t_n,\ \overline t_{n+1}\}\to\{\delta_{k,n}, \delta_{k,n+1}\}^\infty_{k=1}\to\{\gamma_n, \gamma_{n+1}\}\\ \cos\theta(\overline t_n)=0,\ f_k(\delta_{k,n})=0,\ Z(\gamma_n)=0
\end{gathered}
\end{equation}
where
\begin{equation}
\notag
\begin{gathered}
Z(t)=\zeta(1/2+it)e^{i\theta(t)}\\
e^{i\theta(t)}=\pi^{-it/2}\frac{\Gamma(1/4+it/2)}{|\Gamma(1/4+it/2)|}\\
\zeta(s)=\sum^\infty_{n=1} n^{-s}\\
s=\sigma+it
\end{gathered}
\end{equation}
\par 
ii) all zeros of the Hardy function obtained by transformation are real;
\par 
iii) the Hardy function has no other zeros except those obtained by transformation.
\end{theorem}
\par 
In the course of the proof, several expressions were obtained that had not previously been considered by other researchers of the Zeta function.
\par 
Let us briefly review the main stages of the proof.
\par 
First, we have shown that the value of the Zeta function is located in the center of an almost regular polygon, which is formed by partial sums of Dirichlet series, and we have found an expression for the error in calculating the center of a regular polygon through the ratio of the radius of the inlined and outlined circle

\begin{equation}
\notag
\begin{gathered}
\epsilon=r^{(\nu_m)}=r^{(0)}\big(cos(\frac{\pi}{k_m})\big)^{\nu_m}\\
r^{(0)}= \frac{m^{-\sigma}}{2\tan(\frac{\pi}{k_m})}
\end{gathered}
\end{equation}
\par 
Then we have constructed an analytical continuation of the Zeta function, suitable for transformation by a binomial series
\begin{equation}
\notag
\begin{gathered}
\forall\epsilon>0,\ \exists N: \forall k_m>N,\  \nu_m>k_m^{2.15}\\
\Big|\zeta(s)-z^{(\nu_m)}_m\Big|<\epsilon\\
z^{(\nu_m)}_m=\sum^{m+k_m}_{n=1}a_n\delta^{(\nu_m)}_{\lambda_{\omega},n}\\
a_n=n^{-\sigma} \{\cos(t\log(n))-i\sin(t\log(n))\}\\
m=\Big\lfloor\frac{k_m|t|}{2\pi}\Big\rfloor,\ k_m=3,4,5...
\end{gathered}
\end{equation}
\par 
where $\delta^{(\nu_m)}_{\lambda_{\omega},n}$ - coefficients averaging the initial terms of the Dirichlet series $\{a_{\lambda}\}^{m+k_m}_{\lambda=m+1},\ \forall\lambda_{\omega}\in(m+1,\ m+k_m)$.
\begin{equation}
\notag
\begin{gathered}
\delta^{(\nu_m)}_{\lambda_{\omega},n}=
\begin{cases}
n\le m,\ &1\\
m+k_m>n>m &\frac{1}{2^{\nu_m}}\sum\binom {\nu_m} \mu,\ n-m\le (\lambda_{\omega}+\mu)\mod k_m\\
n=m+k_m &\frac{1}{2^{\nu_m}}\sum\binom {\nu_m} \mu,\ (\lambda_{\omega}+\mu)\mod k_m = 0 
\end{cases}\\
\mu=0,1,...\nu_{m}
\end{gathered}
\end{equation}
\par 
We have defined the Generalized Hardy function $Z_{\alpha_\nu}(t)=Re\ \zeta(\alpha_\nu+it)e^{i\theta(t)}$ as the real part of the Hardy function (\ref{hardy_function_1}) on any line $\alpha_\nu+it$ parallel to the critical line $1/2+it$.
\par 
Then used the obtained expression of the analytical continuation of the Zeta function to summarize the Hardy function on any line $\alpha_\nu+it$ parallel to the critical line $1/2+it$
\begin{equation}
\notag
\begin{gathered}
\zeta(\alpha_\nu+it) e^{i\theta(t)}=\lim_{m\to\infty}\overline z^{(\nu_m)}_{\alpha_\nu,m}\\
\overline z^{(\nu_m)}_{\alpha_\nu,m}=\sum^{m+k_m}_{n=1}\overline a_n\delta^{(\nu_m)}_{\lambda_{\omega},n}+i\sum^{m+k_m}_{n=1}\overline b_n\delta^{(\nu_m)}_{\lambda_{\omega},n}\\
\overline a_n=n^{-\alpha_\nu} \cos(\theta(t)-t\log(n)),\ \overline b_n=n^{-\alpha_\nu} \sin(\theta(t)-t\log(n))\\
e^{i\theta(t)}=\pi^{-it/2}\frac{\Gamma(1/4+it/2)}{|\Gamma(1/4+it/2)|}
\end{gathered}\end{equation}
\par 
We used the property of a function of a complex variable. If a limit has a sequence of functions of a complex variable, then the limits also have separate real and imaginary
parts of this sequence.
\begin{equation}
\notag
\begin{gathered}
\forall\epsilon>0,\ \exists N: \forall k_m>N,\  \nu_m>k_m^{2.15}\\  
\Big|Re\ \zeta(\alpha_\nu+it) e^{i\theta(t)}-Re\ \overline z^{(\nu_m)}_{\alpha_\nu,m}\Big|<\epsilon\\
\text{or}\\
\Big|Z_{\alpha_\nu}-z^{(\nu_m)}_{\alpha_\nu,m}\Big|<\epsilon\\
z^{(\nu_m)}_{\alpha_\nu,m}=Re\ \overline z^{(\nu_m)}_{\alpha_\nu,m}=\sum^{m+k_m}_{n=1}\overline a_n\delta^{(\nu_m)}_{\lambda_{\omega},n}\\
m= \Big\lfloor\frac{k_m|t|}{2\pi}\Big\rfloor,\ k_m=3,4,5...
\end{gathered}
\end{equation}
\par 
to summarize the Generalized Hardy function
\begin{equation}
\notag
\begin{gathered}
Z_{\alpha_\nu}=\cos\theta(t)+\lim_{m\to\infty}\sum^{m+k_m}_{n=2}n^{-\alpha_\nu}T_n\\
T_n=\cos(\theta(t)-t\log(n))\delta^{(\nu_m)}_{\lambda_{\omega},n}
\end{gathered}
\end{equation}
\par 
and substituted the expression of the binomial series into the expression for the Generalized Hardy function
\begin{equation}
\notag
n^{-\alpha_\nu}=1+\sum^{\infty}_{k=1}(-1)^k\frac{\alpha_\nu(\alpha_\nu-1)(\alpha_\nu-2)...(\alpha_\nu-k+1)}{k!}\Big(\frac{n-1}{n}\Big)^k
\end{equation}

\par 
We reversed the order of summation under the sign of the limit\footnote{because this is the sequence of functions $\{z^{(\nu_m)}_{\alpha_\nu,m}\}$ converges uniformly, and the binomial series for $|x|<1$ converges absolutely} and got an expression
\begin{equation}
\notag
\begin{gathered}
Z_{\alpha_\nu} =Z_0+\sum^{\infty}_{k=1}(-1)^k \frac{\alpha_\nu(\alpha_\nu-1)(\alpha_\nu-2)...(\alpha_\nu-k+1)}{k!}Q_k\\
Z_0=\cos\theta(t)+\lim_{m\to\infty}\sum^{m+k_m}_{n=2}T_n\\
Q_k=\lim_{m\to\infty}\sum^{m+k_m}_{n=2}\Big(\frac{n-1}{n}\Big)^kT_n,\ k=0,1,2...\\
m=\Big\lfloor\frac{k_m|t|}{2\pi}\Big\rfloor,\ k_m=3,4,5...
\end{gathered}
\end{equation}
\par 
Then, in a few steps, we determined the transformation of the zeros of the Hardy function:
\par 
1) it is easy to see that 
$$\lim_{\alpha_\nu\to\infty}Z_{\alpha_\nu}=\cos\theta(t)$$
\par 
Notice that
$$\theta(t)=\frac{t}{2}\log\frac{t}{2\pi}-\frac{t}{2}-\frac{\pi}{8}+O(1/t)$$
\par 
distance between neighboring zeros of $\cos\theta(t)$ decreases with increasing $t$ therefor we split the numeric axis into a countable number of intervals $\Delta T_\lambda=(t_\lambda, t_{\lambda+1}]$ 
$$t_\lambda=2\pi\lambda^2,\ \lambda=1,\ 2,\ 3\ ...$$
\par 
and define constant $A_\lambda$ for each interval
\begin{equation}
\notag
\begin{gathered}
\forall \Delta T_\lambda=(t_\lambda, t_{\lambda+1}],\ 
t_\lambda=2\pi\lambda^2,\ \lambda=1,\ 2,\ 3\ ...\\
\exists A_\lambda:\forall \hat\alpha_\lambda>A_\lambda\\
|\cos\theta(t) -Z_{\hat\alpha_\lambda}(t)|<\epsilon(A_\lambda),\ t\in \Delta T_\lambda
\end{gathered}
\end{equation}
\par 
2) it is easy to obtain a sequence of functions to transform the Generalized Hardy function 
$$Z_{\alpha^{(\lambda)}_{\nu+1}}\to \{Z_{\alpha^{(\lambda)}_\nu, k}\}^\infty_{k=1}\to Z_{\alpha^{(\lambda)}_{\nu}}$$
between any two lines $\alpha^{(\lambda)}_{\nu}+it$ and $\alpha^{(\lambda)}_{\nu+1}+it$ parallel to the critical line $1/2+it$ for each interval\footnote{$\Delta_\lambda$ denotes the tolerance for the shift zeros of $Z_{\alpha^{(\lambda)}_\nu}(t)$ relative to zeros of $\cos\theta(t)$} $\Delta T_\lambda\pm\Delta_\lambda$
\begin{equation}
\notag
\begin{gathered}
Z_{\alpha^{(\lambda)}_\nu, 0}=Z_{\alpha^{(\lambda)}_{\nu+1}}\\
Z_{\alpha^{(\lambda)}_\nu, k}=Z_{\alpha^{(\lambda)}_\nu, k-1}-P_{\alpha^{(\lambda)}_\nu,k}
\end{gathered}
\end{equation}
\begin{equation}
\notag
\begin{split}
P_{\alpha^{(\lambda)}_\nu,k}=\frac{(-1)^k}{k!}\big\{\alpha^{(\lambda)}_{\nu+1}(\alpha^{(\lambda)}_{\nu+1}-1)(\alpha^{(\lambda)}_{\nu+1}-2)&...(\alpha^{(\lambda)}_{\nu+1}-k+1)\\
-\alpha^{(\lambda)}_{\nu}(\alpha^{(\lambda)}_{\nu}-1)(\alpha^{(\lambda)}_{\nu}-2)&...(\alpha^{(\lambda)}_{\nu}-k+1)\big\}Q_k
\end{split}
\end{equation}
\par 
3) to have possibility to decrease value $\max|P_{\alpha^{(\lambda)}_\nu,k}|$ between any neighboring zeros of $Z_{\alpha^{(\lambda)}_\nu, k-1}$ we fix the $\mu_\lambda$ values 
$$\alpha^{(\lambda)}_1<\alpha^{(\lambda)}_2<\alpha^{(\lambda)}_3<...<\alpha^{(\lambda)}_{\nu}<...<\alpha^{(\lambda)}_{\mu_\lambda}$$ 
between the values $\sigma=1/2$ and $\sigma=\hat\alpha_\lambda>A_\lambda$ for each interval $\Delta T_\lambda\pm\Delta_\lambda$ and define a sequence of functions

\begin{equation}
\notag
\begin{split}
\{f^{(\lambda)}_k\}^\infty_{k=1} &=\big\{\cos\theta\to Z_{\hat\alpha_\lambda}\to \{Z_{\alpha^{(\lambda)}_{\mu_\lambda}, k}\}^\infty_{k=1}\\
&\to Z_{\alpha^{(\lambda)}_{\mu_\lambda}}\to \{Z_{\alpha^{(\lambda)}_{\mu_\lambda-1}, k}\}^\infty_{k=1}\\
...&\to Z_{\alpha^{(\lambda)}_{\nu+1}}\to \{Z_{\alpha^{(\lambda)}_\nu, k}\}^\infty_{k=1}\to Z_{\alpha^{(\lambda)}_{\nu}}...\\
&\to Z_{\alpha^{(\lambda)}_3}\to  \{Z_{\alpha^{(\lambda)}_2, k}\}^\infty_{k=1}\\
&\to Z_{\alpha^{(\lambda)}_2}\to  \{Z^{(\lambda)}_{\alpha^{(\lambda)}_1, k}\}^\infty_{k=1}\to Z\big\}
\end{split}
\end{equation}
\par 
as result we define sequence of transformation of zeros $\{f^{(\lambda)}_k\}^\infty_{k=1}$
\begin{equation}
\notag
\begin{split}
&\{\overline t^{(\lambda)}_n,\ \overline t^{(\lambda)}_{n+1}\}\to \{\delta_{\hat\alpha_\lambda,n},\ \delta_{\hat\alpha_\lambda,n+1}\}\to \{\delta_{\alpha^{(\lambda)}_{\mu_\lambda},k,n},\ \delta_{\alpha^{(\lambda)}_{\mu_\lambda},k,n+1}\}^\infty_{k=1}\\
&\to \{\delta_{\alpha^{(\lambda)}_{\mu_\lambda},n},\ \delta_{\alpha^{(\lambda)}_{\mu_\lambda},n+1}\}\to \{\delta_{\alpha^{(\lambda)}_{\mu_\lambda-1},k,n},\ \delta
_{\alpha^{(\lambda)}_{\mu_\lambda-1},k,n+1}\}^\infty_{k=1}\\
...&\to \{\delta_{\alpha^{(\lambda)}_{\nu+1},n},\ \delta_{\alpha^{(\lambda)}_{\nu+1},n+1}\}\to \{\delta_{\alpha^{(\lambda)}_\nu,k,n},\ \delta_{\alpha^{(\lambda)}_\nu,k,n+1}\}^\infty_{k=1}\\
&\to \{\delta_{\alpha^{(\lambda)}_\nu,n},\ \delta_{\alpha^{(\lambda)}_\nu,n+1}\}...\\
&\to \{\delta_{\alpha^{(\lambda)}_3,n},\ \delta_{\alpha^{(\lambda)}_3,n+1}\}\to \{\delta_{\alpha^{(\lambda)}_2,k,n},\ \delta_{\alpha^{(\lambda)}_2,k,n+1}\}^\infty_{k=1}\\
&\to \{\delta_{\alpha^{(\lambda)}_2,n},\ \delta_{\alpha^{(\lambda)}_2,n+1}\}\to \{\delta^{(\lambda)}_{\alpha_1,k,n},\ \delta^{(\lambda)}_{\alpha_1,k,n+1}\}^\infty_{k=1}\to \{\gamma^{(\lambda)}_n,\ \gamma^{(\lambda)}_{n+1}\}
\end{split}
\end{equation}
\par 
where
\begin{equation}
\notag
\begin{gathered}
\cos\theta(\overline t^{(\lambda)}_n)=0,\ Z_{\hat\alpha_\lambda}(\delta_{\hat\alpha_\lambda,n})=0,\ Z_{\alpha^{(\lambda)}_{\nu}}(\delta_{\alpha^{(\lambda)}_{\nu},n})=0\\
Z_{\alpha^{(\lambda)}_\nu, k}(\delta_{\alpha^{(\lambda)}_\nu,k,n})=0,\ Z(\gamma^{(\lambda)}_n)=0
\end{gathered}
\end{equation}
\par 
thus, we proved the first part of the theorem on the transformation of zeros $\cos\theta(t)$ into zeros of the Hardy function.
\par 
To prove the second part of the theorem, we reviewed different options of ranges between points where the function $Z_{\alpha^{(\lambda)}_\nu, k-1}$ changes sign and $P_{\alpha^{(\lambda)}_\nu,k}>0$ or $P_{\alpha^{(\lambda)}_\nu,k}<0$.
\par 
It is four options for ranges $\{1,\ 2,\ 3,\ 4\}$, where $Z_{\alpha^{(\lambda)}_\nu, k-1}>0$ and the same four options for ranges $\{1',\ 2',\ 3',\ 4'\}$, where $Z_{\alpha^{(\lambda)}_\nu, k-1}<0$ that can change each other in a restricted number of ways
\begin{equation}
\notag
\begin{gathered}
\xymatrix{
&&(1)&&(3)\\
(1)\ar[r]^A\ar[dr]&(1')\ar[r]^B\ar[ur]&(3)\ar[r]^B\ar[dr]&(4')\ar[r]^B\ar[ur]&(1)\\
(4)\ar[r]^A\ar[ur]&(3')\ar[r]^A\ar[dr]&(2)\ar[r]^B\ar[ur]&(2')\ar[r]^A\ar[dr]&(4)\\
&&(4)&&(2)
}
\end{gathered}
\end{equation}
\par 
We have selected the following values $\Delta\alpha^{(\lambda)}_\nu=\alpha^{(\lambda)}_{\nu+1}-\alpha^{(\lambda)}_\nu$ so that the condition $\max|Z_{\alpha^{(\lambda)}_{\nu,k-1}}|>\max|P_{\alpha^{(\lambda)}_\nu,k}|$ is met on each range between any neighboring zeros of $Z_{\alpha^{(\lambda)}_\nu, k-1}$.
\par 
Then we got two relation options
\par 
A) $\delta_{\alpha^{(\lambda)}_\nu,k,n}<\delta_{\alpha^{(\lambda)}_\nu,k-1,n}$
\par 
B) $\delta_{\alpha^{(\lambda)}_\nu,k,n}>\delta_{\alpha^{(\lambda)}_\nu,k-1,n}$
\par 
which correspond to a strictly decreasing
$$\delta_{\alpha^{(\lambda)}_\nu,k,n}>\delta_{\alpha^{(\lambda)}_\nu,k+1,n}>\delta_{\alpha^{(\lambda)}_\nu,k+2,n}>...$$
\par 
or a strictly increasing sequence of zeros 
$$\delta_{\alpha^{(\lambda)}_\nu,k,n}<\delta_{\alpha^{(\lambda)}_\nu,k+1,n}<\delta_{\alpha^{(\lambda)}_\nu,k+2,n}<...$$
\par 
We have shown that all zeros of the sequence are real, then the real zeros of the Generalized Hardy function and, by induction, the zeros of the Hardy function obtained by transforming the zeros of the function are also real.
\par 
The proof of the third part is based on 
$$\theta(t)=\frac{t}{2}\log\frac{t}{2\pi}-\frac{t}{2}-\frac{\pi}{8}+O(1/t)$$
we can determine the exact number of zeros $\cos\theta(t)$ for any interval $(0,T]$
$$N_\theta(T)=\frac{T}{2\pi}\log\frac{T}{2\pi}-\frac{T}{2\pi}-\frac{1}{8}+O(1/T)$$
\par 
what corresponds to the Riemann-Mangoldt formula for the number of zeros of the Zeta function in the same interval 
$$N(T)=\frac{T}{2\pi}\log\frac{T}{2\pi}-\frac{T}{2\pi}+O(\log T)$$
\par 
Now let us review the stages of the proof in detail.
\par 
\section{Summation of Zeta function}
The Zeta function is defined by the Dirichlet series
\begin{equation}
\notag
\begin{split}
\zeta(s)&=\sum^\infty_{n=1} n^{-s}\\
&=\sum^\infty_{n=1} n^{-\sigma-it}\\
&=\sum^\infty_{n=1}n^{-\sigma } \{\cos(t\log(n))-i\sin(t\log(n))\}
\end{split}
\end{equation}
\par 
Obviously, the Dirichlet series diverges at $\sigma\le 1$.
\par 
Using the functional equation (\ref{zeta_func_eq}), the values  of the Zeta function at $\sigma<0$ can be calculated by knowing its values at $\sigma>1$, whatever leaving a narrow strip of $0\le\sigma\le 1$ which is called critical where other methods for calculating Zeta function values are required.
\par 
First, we shown that the value of the Zeta function is located in the center of an almost regular polygon, which is formed by partial sums of Dirichlet series, and we found an expression for the error in calculating the center of a regular polygon through the ratio of the radius of the inlined and outlined circle.
\par 
Then we will obtain an analytical continuation of $\zeta(s)$ in a form convenient for transformation by a binomial series.
\par 
We define
\begin{equation}\label{zeta_2}a_n=n^{-\sigma } \{\cos(t\log(n))-i\sin(t\log(n))\}\end{equation}
\par 
and review the sequence of partial sums $\{s_m\}$ of the Dirichlet series on the complex plane
\begin{equation}
\notag
s_m=\sum^m_{n=1}a_n
\end{equation}
\par 
We can see (Fig. \ref{fig:polygon_1}) that when
\begin{equation}
\notag
m= \Big\lfloor\frac{kt}{2\pi}\Big\rfloor,\ k=3,4,5...\end{equation}
\par 
the segments between the points corresponding to the k partial sums form an almost\footnote{the relative angle between the segments $a_n$ that connect the points corresponding to the partial sums $s_m$ of the Dirichlet series increases slowly, while $|a_n|$ decreases slowly}  regular $k$ polygon.
\begin{center}
\includegraphics[scale=0.5]{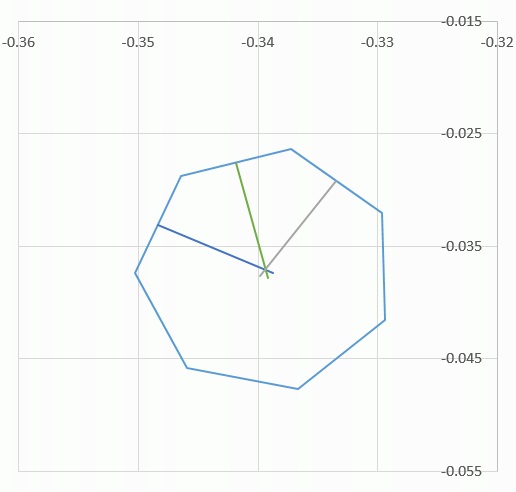}
\end{center}
\begin{figure}[ht!]
\centering
\caption{Polygon}
\label{fig:polygon_1}
\end{figure}
\par 
Obviously, according to the asymptotic expression (\ref{zeta_5}), the point corresponding to the value of the Zeta function \cite{TI} for sufficiently large $N$ is located in the center of this polygon
\begin{equation}
\label{zeta_5}
\begin{gathered}
\zeta(s)\sim\sum^N_{n=1}n^{-s}-\frac{N^{1-s}}{1-s}-\frac{N^{-s}}{2}\\
|s|>-1\\
s=\sigma+it
\end{gathered}
\end{equation}
\par 
Indeed, the last term (\ref{zeta_5}) returns us to the middle of the last segment, and the second term is geometrically perpendicular to the last term because the numerator of the second term in exponential form has the same argument as the third term
\begin{equation}
\notag
\begin{split}
N^{1-s} & =N^{1-\sigma}e^{ -it\log(N)}\\
N^{-s} & =N^{-\sigma}e^{ -it\log(N)}
\end{split}
\end{equation}
\par 
and the denominator of the second term (\ref{zeta_5}) is almost perpendicular to the numerator in accordance with the denominator argument in exponential form
$$\frac{1}{1-s}=\frac{1}{\sqrt{(1-\sigma)^2+t^2)}}e^{ i\arctan\frac{t}{1-\sigma}}$$
\par 
because
\par 
\begin{equation}
\notag
\begin{split}
\arctan\frac{t}{1-\sigma} & \to \frac{\pi}{2}\\
t & \to\infty
\end{split}
\end{equation}
\par 
It is obvious that
\begin{equation}
\notag
\begin{split}
\frac{N^{1-\sigma}}{\sqrt{(1-\sigma)^2+t^2)}} & \to\frac{N^{1-\sigma}}{t}\\
t & \to\infty
\end{split}
\end{equation}
\par 
Since for the last term (\ref{zeta_5}) is true 
$$N= \Big\lfloor\frac{kt}{2\pi}\Big\rfloor\to \frac{N}{t}=\frac{k}{2\pi}$$
\par 
And also that
\begin{equation}
\notag
\begin{split}
\tan\Big(\frac{\pi}{k}\Big) & \sim\frac{\pi}{k}\\
k & \to\infty
\end{split}
\end{equation}
\par 
we get confirmation that the point corresponding to the value of the Zeta function is located in the center of the k polygon formed by the segments between the points corresponding to the k partial sums in geometric form
\begin{equation}
\notag
\begin{split}
\frac{N^{1-\sigma}}{\sqrt{(1-\sigma)^2+t^2)}} & \to\frac{N^{-\sigma}}{2\tan\big(\frac{\pi}{k}\big)}\\
t & \to\infty
\end{split}
\end{equation}
\par 
since the resulting expression corresponds to the radius of a circle inlined in a regular k polygon
\begin{equation}
\label{inscribed_circle}
r=\frac{a}{2\tan(\frac{\pi}{k})}
\end{equation}
\par
\begin{figure}[ht!]
\centering
\includegraphics[scale=0.5]{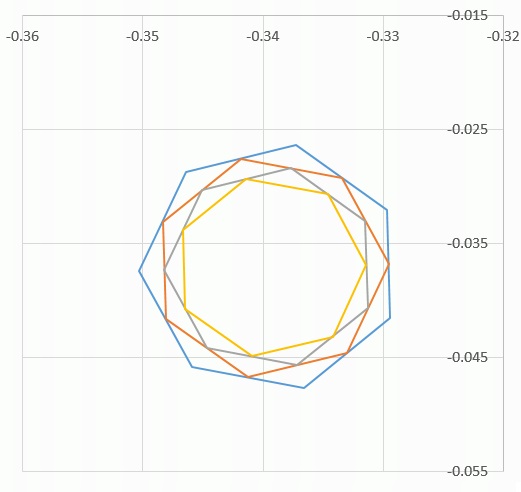}
\caption{Inlined polygons}
\label{fig:polygon_2}
\end{figure}

 
Calculating $\nu_m$ times the average values of the partial sums  $\{s^{(\nu_m)}_{\lambda}\}^{m+k_m}_{\lambda=m+1
}$ that correspond to $k_m$  top points of the k polygon, we obtain a sequence of inlined k polygons (Fig. \ref{fig:polygon_2}), which obviously converges to the center of the original k polygon
\par 
Since for an outer $k$ polygon, the circle is inlined, and for an inlined $k$ polygon, the same circle is outlined, then
$$r^{(0)}=R^{(1)}=\frac{a^{(1)}}{2\sin(\frac{\pi}{k_m})}$$
\par 
using
$$r^{(1)}=\frac{a^{(1)}}{2\tan(\frac{\pi}{k_m})}$$
\par 
we obtain the recurrence relation
$$r^{(1)}=r^{(0)}\cos(\frac{\pi}{k_m})$$
\par 
and the expression for the error of calculating the center of the k polygon by averaging partial sums \footnote{the coefficient 2.15 is obtained by calculating}
\begin{equation}
\label{epsilon_1}
\begin{split}
\epsilon=r^{(\nu_m)} & =r^{(0)}\big(cos(\frac{\pi}{k_m})\big)^{\nu_m}\\
\nu_m&>k_m^{2.15} \\
r^{(0)} & = \frac{m^{-\sigma}}{2\tan(\frac{\pi}{k_m})}\\
m & = \Big\lfloor\frac{k_mt}{2\pi}\Big\rfloor,\ k_m=3,4,5...
\end{split}
\end{equation}
\par 
As well as expressions\footnote{obviously any top k of the polygon $s^{(\nu_m)}_{\lambda_{\omega}}$ will be located in a circle with radius $\epsilon$} to calculate the coefficients of $\delta^{(\nu_m)}_{\lambda_{\omega},n}$ averaging the initial terms of the Dirichlet series $\{a_{\lambda}\}^{m+k_m}_{\lambda=m+1}$
\begin{equation}
\notag
\begin{split}
s^{(0)}_{\lambda_{\omega}} & =s_{\lambda_{\omega}} \\ 
s^{(1)}_{\lambda_{\omega}} & =\frac{s_{\lambda_{\omega}} + s_{\lambda_{\omega}+1}}{2} \\ 
s^{(2)}_{\lambda_{\omega}} & =\frac{ s^{(1)}_{\lambda_{\omega}} + s^{(1)}_{\lambda_{\omega}+1}}{2}\\
& ...\\
s^{(\nu_m)}_{\lambda_{\omega}} & =\frac{ s^{(\nu_m-1)}_{\lambda_{\omega}} + s^{(\nu_m-1)}_{\lambda_{\omega}+1}}{2}
\end{split}
\end{equation}
\par 
or through partial sums
\begin{equation}
\notag
\begin{split}
s^{(2)}_{\lambda_{\omega}} & =\frac{ s_{\lambda_{\omega}} + 2s_{\lambda_{\omega}+1}+s_{\lambda_{\omega}+2}}{4}\\
s^{(3)}_{\lambda_{\omega}} & =\frac{ s_{\lambda_{\omega}} + 3s_{\lambda_{\omega}+1}+3s_{\lambda_{\omega}+2}+s_{\lambda_{\omega}+3}}{8}\\
& ...\\
s^{(\nu_m)}_{\lambda_{\omega}} & =\frac{ s_{\lambda_{\omega}} + \binom {\nu_m} 1 s_{\lambda_{\omega}+1}+\binom {\nu_m} 2 s_{\lambda_{\omega}+2}+...+\binom {\nu_m} {\mu} s_{\lambda_{\omega}+\mu}+...+s_{\lambda_{\omega}+\nu_m}}{2^{\nu_m}}
\end{split}
\end{equation}
\par 
then, using the conditional \footnote{since $k_m<\nu_m$, the partial sums of $s_{\lambda_{\omega}+\mu}$ form a cyclic sequence that begins with $s_{\lambda_{\omega}}$ then each time after $s_{m+k_m}$ continues with $s_{m+1}$ until it reaches $s_{\lambda_{\omega}+\nu_m}$, so we sum the binomial coefficients if a term of the original series $a_n$ sums to $s_{\lambda_{\omega}+\mu}$, because we chose an arbitrary vertex $s^{(\nu_m)}_{\lambda_{\omega}}$} numbering $s_{\lambda_{\omega}+\mu}$ we can calculate the coefficients to sum the terms of the original series
\begin{equation}
\label{delta_nu}
\begin{gathered}
\delta^{(\nu_m)}_{\lambda_{\omega},n}=
\begin{cases}
n\le m,\ &1\\
m+k_m>n>m &\frac{1}{2^{\nu_m}}\sum\binom {\nu_m} \mu,\ n-m\le (\lambda_{\omega}+\mu)\mod k_m\\
n=m+k_m &\frac{1}{2^{\nu_m}}\sum\binom {\nu_m} \mu,\ (\lambda_{\omega}+\mu)\mod k_m = 0 
\end{cases}\\
\mu=0,1,...\nu_{m}
\end{gathered}
\end{equation}
\par 
Define the sequence of functions $\{z^{(\nu_m)}_m\}$, where $a_n$ and $\delta^{(\nu_m)}_{\lambda_{\omega},n}$ are defined by (\ref{zeta_2}) and (\ref{delta_nu}), respectively
\begin{equation}
\notag
z^{(\nu_m)}_m=\sum^{m+k_m}_{n=1}a_n\delta^{(\nu_m)}_{\lambda_{\omega},n}
\end{equation}
\par 
Obviously, based on (\ref{zeta_5}, \ref{epsilon_1}, \ref{delta_nu}), that the sequence of functions $\{z^{(\nu_m)}_m\}$ converges to the value of the Zeta function
\begin{equation}
\label{epsilon_2}
\begin{gathered}
\forall\epsilon>0,\ \exists N: \forall k_m>N,\  \nu_m>k_m^{2.15}\\
\Big|\zeta(s)-z^{(\nu_m)}_m\Big|<\epsilon\\
m=\Big\lfloor\frac{k_m|t|}{2\pi}\Big\rfloor,\ k_m=3,4,5...
\end{gathered}
\end{equation}
\par 
in addition, the sequence of functions $\{z^{(\nu_m)}_m\}$ converges uniformly on the entire complex plane, except for a narrow strip along the real axis
\begin{equation}
\notag
\begin{gathered}
\forall\epsilon>0,\ |t|>|t_0|>0,\ \sigma\subset\mathbb{R}\\
\exists N: \forall k_p,k_q>N,\  \nu_p>k_p^{2.15},\ \nu_q>k_q^{2.15}\\
\Big|z^{(\nu_p)}_p-z^{(\nu_q)}_q\Big|<\epsilon\\ 
p=\Big\lfloor\frac{k_p|t|}{2\pi}\Big\rfloor,\ q=\Big\lfloor\frac{k_q|t|}{2\pi}\Big\rfloor
\end{gathered}
\end{equation}
\par 
\section{Generalized Hardy function}
\par 
Review the Hardy function (\ref{hardy_function_1}) on any line $\alpha_\nu+it$ parallel to the critical line $1/2+it$
\begin{equation}
\notag
\begin{split}
\zeta(\alpha_\nu+it) e^{i\theta(t)} & =\sum^\infty_{n=1}n^{-\alpha_\nu-it} e^{i\theta(t)}\\
& =\sum^\infty_{n=1}n^{-\alpha_\nu } e^{i(\theta(t)-t\log(n))}\\
& =\sum^\infty_{n=1}n^{-\alpha_\nu } \cos(\theta(t)-t\log(n))+i\sum^\infty_{n=1}n^{-\alpha_\nu }\sin(\theta(t)-t\log(n))
\end{split}
\end{equation}
\par 
\begin{definition}
Let's define $Z_{\alpha_\nu}$ by the Generalized Hardy function only the real part of the Hardy function on any line $\alpha_\nu+it$ parallel to the critical line $1/2+it$, because on the critical line the imaginary part of the Hardy function is zero
\begin{equation}
\label{generalized_hardy_function_1}
\begin{gathered}
Z_{\alpha_\nu}=Re\ \zeta(\alpha_\nu+it)e^{i\theta(t)}=\sum^\infty_{n=1}n^{-\alpha_\nu }\cos(\theta(t)-t\log(n))\\
e^{i\theta(t)}=\pi^{-it/2}\frac{\Gamma(1/4+it/2)}{|\Gamma(1/4+it/2)|}\\
\zeta(s)=\sum^\infty_{n=1} n^{-s}\\
s=\alpha_\nu+it
\end{gathered}
\end{equation}
\end{definition}
\par 
Obviously, we are taking the projection $u(1/2+\Delta+it)$ of the Zeta function onto the real axis of the coordinate system of the complex plane, which we have rotated by the angle $\theta(t)$
\begin{center}
\includegraphics[scale=0.7]{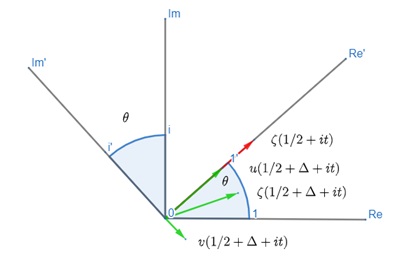}
\end{center}
\begin{figure}[ht!]
\centering
\caption{Generalized Hardy function $u(1/2+\Delta+it)$}
\label{fig:rotation_2}
\end{figure}
\par 
In according to the functional equation (\ref{zeta_func_eq}) such a function always has a pair $u(1/2-\Delta+it)$
\begin{center}
\includegraphics[scale=0.7]{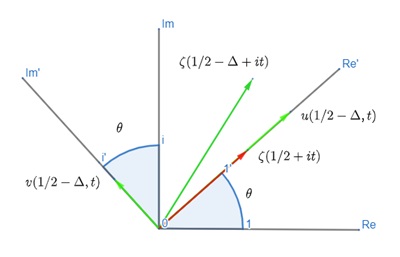}
\end{center}
\begin{figure}[ht!]
\centering
\caption{Paired Generalized Hardy function $u(1/2-\Delta+it)$}
\label{fig:rotation_3}
\end{figure}
\par 
\section{Summation of Generalized Hardy function}
We define the Hardy function on any line $\alpha_\nu+it$ parallel to the critical line $1/2+it$ through the sum of the Zeta function (\ref{epsilon_2}), where $\delta^{(\nu_m)}_{\lambda_{\omega},n}$ are still being defined (\ref{delta_nu})
\begin{equation}
\notag
\begin{gathered}
\zeta(\alpha_\nu+it) e^{i\theta(t)}=\lim_{m\to\infty}\overline z^{(\nu_m)}_{\alpha_\nu,m}\\
\overline z^{(\nu_m)}_{\alpha_\nu,m}=\sum^{m+k_m}_{n=1}\overline a_n\delta^{(\nu_m)}_{\lambda_{\omega},n}+i\sum^{m+k_m}_{n=1}\overline b_n\delta^{(\nu_m)}_{\lambda_{\omega},n}\\
\overline a_n=n^{-\alpha_\nu } \cos(\theta(t)-t\log(n)),\ \overline b_n=n^{-\alpha_\nu } \sin(\theta(t)-t\log(n))\\
e^{i\theta(t)}=\pi^{-it/2}\frac{\Gamma(1/4+it/2)}{|\Gamma(1/4+it/2)|}
\end{gathered}
\end{equation}
Since
$$|e^{i\theta(t)}|=1$$
\par 
obviously, what is the sequence of functions $\{\overline z^{(\nu_m)}_{\alpha_\nu,m}\}$ converges to the value of the Hardy function on any line $\alpha_\nu+it$ parallel to the critical line $1/2+it$
\begin{equation}
\label{epsilon_4}
\begin{gathered}
\forall\epsilon>0,\ \exists N: \forall k_m>N,\  \nu_m>k_m^{2.15}\\  \Big|\zeta(\alpha_\nu+it) e^{i\theta(t)}-\overline z^{(\nu_m)}_{\alpha_\nu,m}\Big|<\epsilon\\
m= \Big\lfloor\frac{k_m|t|}{2\pi}\Big\rfloor,\ k_m=3,4,5...
\end{gathered}
\end{equation}
\par 
in addition, the sequence of functions $\{\overline z^{(\nu_m)}_m\}$ converges uniformly on the entire complex plane, except for a narrow strip along the real axis
\begin{equation}
\notag
\begin{gathered}
\forall\epsilon>0,\ |t|>|t_0|>0,\ \alpha_\nu\subset\mathbb{R}\\
\exists N: \forall k_p,k_q>N,\  \nu_p>k_p^{2.15},\ \nu_q>k_q^{2.15}\\  \Big|\overline z^{(\nu_p)}_{\alpha_\nu,p}-\overline z^{(\nu_q)}_{\alpha_\nu,q}\Big|<\epsilon \\
p=\Big\lfloor\frac{k_p|t|}{2\pi}\Big\rfloor,\ q=\Big\lfloor\frac{k_q|t|}{2\pi}\Big\rfloor
\end{gathered}
\end{equation}
\par 
We use the property of the limit of a sequence of functions of a complex variable. If the limit has a sequence of functions of a complex variable, then the limits also have separate real and imaginary
parts of this sequence.
\begin{equation}
\notag
\begin{split}
\zeta(\alpha_\nu+it) e^{i\theta(t)}&=\lim_{m\to\infty}\overline z^{(\nu_m)}_{\alpha_\nu,m}\\&=\lim_{m\to\infty}\sum^{m+k_m}_{n=1}\overline a_n\delta^{(\nu_m)}_{\lambda_{\omega},n}\/+i\lim_{m\to\infty}\sum^{m+k_m}_{n=1}\overline b_n\delta^{(\nu_m)}_{\lambda_{\omega},n}
\end{split}\end{equation}
\par 
We define the sequence of functions $\{z^{(\nu_m)}_{\alpha_\nu,m}\}$ as the real part of a sequence of functions $\{\overline z^{(\nu_m)}_{\alpha_\nu,m}\}$
\begin{equation}
\notag
z^{(\nu_m)}_{\alpha_\nu,m}=Re\ \overline z^{(\nu_m)}_{\alpha_\nu,m}=\sum^{m+k_m}_{n=1}\overline a_n\delta^{(\nu_m)}_{\lambda_{\omega},n}\end{equation}
\par 
Obviously, the sequence of functions $\{z^{(\nu_m)}_{\alpha_\nu,m}\}$ converges to the value of the Generalized Hardy function by definition
\begin{equation}
\notag
\begin{gathered}
\forall\epsilon>0,\ \exists N: \forall k_m>N,\  \nu_m>k_m^{2.15}\\  \Big|Re\ \zeta(\alpha_\nu+it) e^{i\theta(t)}-Re\ \overline z^{(\nu_m)}_{\alpha_\nu,m}\Big|<\epsilon\\
\text{or}\\
\Big|Z_{\alpha_\nu}-z^{(\nu_m)}_{\alpha_\nu,m}\Big|<\epsilon\\
m= \Big\lfloor\frac{k_m|t|}{2\pi}\Big\rfloor,\ k_m=3,4,5...
\end{gathered}
\end{equation}
\par 
in addition, the sequence of functions $\{z^{(\nu_m)}_{\alpha_\nu,m}\}$ converges uniformly on the entire complex plane, except for a narrow strip along the real axis
\begin{equation}
\notag
\begin{gathered}
\forall\epsilon>0,\ |t|>|t_0|>0,\ \alpha_\nu\subset\mathbb{R}\\
\exists N: \forall k_p,k_q>N,\  \nu_p>k_p^{2.15},\ \nu_q>k_q^{2.15}\\
\Big|Re\ \overline z^{(\nu_p)}_{\alpha_\nu,p}-Re\ \overline z^{(\nu_q)}_{\alpha_\nu,q}\Big|<\epsilon \\
\text{or}\\
\Big|z^{(\nu_p)}_{\alpha_\nu,p}- z^{(\nu_q)}_{\alpha_\nu,q}\Big|<\epsilon \\
p=\Big\lfloor\frac{k_p|t|}{2\pi}\Big\rfloor,\  q=\Big\lfloor\frac{k_q|t|}{2\pi}\Big\rfloor
\end{gathered}
\end{equation}
\par 
Let
\begin{equation}\label{zeta_7}T_n=\cos(\theta(t)-t\log(n))\delta^{(\nu_m)}_{\lambda_{\omega},n}\end{equation}
\par 
Then we get the expression for the Generalized Hardy function suitable for transformation by a binomial series
\begin{equation}\label{zeta_8}\begin{gathered}
Z_{\alpha_\nu}=\cos\theta(t)+\lim_{m\to\infty}\sum^{m+k_m}_{n=2}n^{-\alpha_\nu }T_n\\
m=\Big\lfloor\frac{k_m|t|}{2\pi}\Big\rfloor,\ k_m=3,4,5...
\end{gathered}\end{equation}
\par 
\section{Binomial series and Generalized Hardy function}
Let
\begin{equation}\label{binomial_series_1}n^{-\alpha_\nu}=\Big(1-\frac{n-1}{n}\Big)^{\alpha_\nu}\end{equation}
\par 
It is obvious that
\begin{equation}\notag\frac{n-1}{n}<1\end{equation}
\par 
Therefore, we can replace the expression (\ref{binomial_series_1}) with a binomial series
\begin{equation}
\notag
(1-x)^{\alpha_\nu}=\sum^{\infty}_{k=0}(-1)^k\frac{\alpha_\nu(\alpha_\nu-1)(\alpha_\nu-2)...(\alpha_\nu-k+1)}{k!}x^k\end{equation}
\par 
Then
\begin{equation}\label{binomial_series_4}n^{-\alpha_\nu}=1+\sum^{\infty}_{k=1}(-1)^k\frac{\alpha_\nu(\alpha_\nu-1)(\alpha_\nu-2)...(\alpha_\nu-k+1)}{k!}\Big(\frac{n-1}{n}\Big)^k\end{equation}
\par 
We substitute (\ref{binomial_series_4}) in the expression for the Generalized Hardy function (\ref{zeta_8})
\begin{equation}
\notag
Z_{\alpha_\nu}=\cos\theta(t)+\lim_{m\to\infty}\sum^{m+k_m}_{n=2}\Big\{1+\sum^{\infty}_{k=1}(-1)^k\frac{\alpha_\nu(\alpha_\nu-1)(\alpha_\nu-2)...(\alpha_\nu-k+1)}{k!}\Big(\frac{n-1}{n}\Big)^k\Big\}T_n
\end{equation}
\par 
Since the sequence of functions $\{z^{(\nu_m)}_{\alpha_\nu,m}\}$ converges uniformly, and the binomial series for $|x|<1$ converges absolutely, we can change the order of summation under the sign of the limit 
\begin{equation}
\notag
Z_{\alpha_\nu} =\cos\theta(t)+\lim_{m\to\infty}\sum^{m+k_m}_{n=2}T_n+\sum^{\infty}_{k=1}(-1)^k\frac{\alpha_\nu(\alpha_\nu-1)(\alpha_\nu-2)...(\alpha_\nu-k+1)}{k!}\lim_{m\to\infty}\sum^{m+k_m}_{n=2}\Big(\frac{n-1}{n}\Big)^kT_n\end{equation}
\par 
Let
\begin{equation}\label{generalized_hardy_function_4}Q_k=\lim_{m\to\infty}\sum^{m+k_m}_{n=2}\Big(\frac{n-1}{n}\Big)^kT_n,\ k=0,1,2...\end{equation}
\par 
Then the sequence of functions $\{q^{(\nu_m)}_m\}$, where $T_n$ is defined by (\ref{zeta_7})
\begin{equation}
\notag
\begin{gathered}
q^{(\nu_m)}_m=\sum^{m+k_m}_{n=2}\Big(\frac{n-1}{n}\Big)^kT_n\\
m=\Big\lfloor\frac{k_m|t|}{2\pi}\Big\rfloor,\ k_m=3,4,5...
\end{gathered}\end{equation}
\par 
as it follows form (\ref{epsilon_1}, \ref{delta_nu}) and
$$\Big (\frac{n-1}{n}\Big)^k<1$$
\par 
converges to the function $Q_k$ by definition
\begin{equation}
\notag
\begin{gathered}
\forall\epsilon>0,\ \exists N: \forall k_m>N,\  \nu_m>k_m^{2.15}\\  \Big|Q_k-q^{(\nu_m)}_m\Big|<\epsilon\\ m= \Big\lfloor\frac{k_m|t|}{2\pi}\Big\rfloor,\ k_m=3,4,5...
\end{gathered}
\end{equation}
\par 
In addition, it converges uniformly on the entire complex plane, except for a narrow strip along the real axis.
\begin{equation}
\notag
\begin{gathered}
\forall\epsilon>0,\ |t|>|t_0|>0,\ \alpha_\nu\subset\mathbb{R}\\
\exists N: \forall k_p,k_q>N,\  \nu_p>k_p^{2.15},\ \nu_q>k_q^{2.15}\\  
\Big|q^{(\nu_p)}_p-q^{(\nu_q)}_q\Big|<\epsilon \\ 
p= \Big\lfloor\frac{k_p|t|}{2\pi}\Big\rfloor,
q= \Big\lfloor\frac{k_q|t|}{2\pi}\Big\rfloor
\end{gathered}
\end{equation}
\par 
It is obvious that
\begin{equation}\label{generalized_hardy_function_5}Z_0=\cos\theta(t)+\lim_{m\to\infty}\sum^{m+k_m}_{n=2}T_n\end{equation}
Then we get the expression for analyzing the zeros of the Generalized Hardy function
\begin{equation}\label{generalized_hardy_function_6}Z_{\alpha_\nu} =Z_0+\sum^{\infty}_{k=1}(-1)^k \frac{\alpha_\nu(\alpha_\nu-1)(\alpha_\nu-2)...(\alpha_\nu-k+1)}{k!}Q_k\end{equation}
\par 
\section{The transformation of zeros \texorpdfstring{$\cos\theta(t)$}{} into zeros of the Hardy function}
From the first expression (\ref{zeta_8}) for the Generalized Hardy function, it is easy to see that $$\lim_{\alpha_\nu\to\infty}Z_{\alpha_\nu}=\cos\theta(t)$$
\par 
Since \cite{VK, TI}
$$\theta(t)=\frac{t}{2}\log\frac{t}{2\pi}-\frac{t}{2}-\frac{\pi}{8}+O(1/t)$$
distance between neighboring zeros of $\cos\theta(t)$ decreases with increasing $t$, therefor we split the numeric axis into a countable number of intervals $\Delta T_\lambda=(t_\lambda, t_{\lambda+1}]$ 
$$t_\lambda=2\pi\lambda^2,\ \lambda=1,\ 2,\ 3\ ...$$
\par 
We can define constant $A_\lambda$ at each interval $\Delta T_\lambda$ based on the constant relationship 
$$H_\lambda B^{-1}=2^{-A_\lambda}$$
\par 
where
$$\\H_\lambda=\frac{2\pi}{\log\text{\large $\frac{t_{\lambda+1}}{2\pi}$}}$$
then
\begin{equation}
\notag
\begin{gathered}
\forall \Delta T_\lambda=(t_\lambda, t_{\lambda+1}],\ 
t_\lambda=2\pi\lambda^2,\ \lambda=1,\ 2,\ 3\ ...\\
\exists A_\lambda:\forall \hat\alpha_\lambda>A_\lambda\\
|\cos\theta(t) -Z_{\hat\alpha_\lambda}(t)|<\epsilon(A_\lambda),\ t\in \Delta T_\lambda
\end{gathered}
\end{equation}
\par 
where
$$A_\lambda=\frac{\log B-\log2\pi+\log(2\log(\lambda+1))}{\log2}$$
\par 
Since $\sin\theta(\overline t_n)=1$ for $\cos\theta(\overline t_n)=0$ and
\begin{equation}
\notag
\begin{gathered}
|\overline t_n -\delta_{\hat\alpha_\lambda,n}|<\epsilon_{\delta}\\
\epsilon_{\delta}=\epsilon(A_\lambda)\sin\theta(\overline t_n)
\end{gathered}
\end{equation}
\par 
Then all zeros of $Z_{\hat\alpha_\lambda}(t)$ are real, because all zeros of $\cos\theta(t)$ are real.
\par 
From the second expression (\ref{generalized_hardy_function_6})  it is easy to obtain a sequence of functions $\{Z_{\alpha^{(\lambda)}_\nu,k}\}$ to transform the Generalized Hardy function
$$Z_{\alpha^{(\lambda)}_{\nu+1}}\to \{Z_{\alpha^{(\lambda)}_\nu, k}\}^\infty_{k=1}\to Z_{\alpha^{(\lambda)}_{\nu}}$$
between any two lines $\alpha^{(\lambda)}_{\nu}+it$ and $\alpha^{(\lambda)}_{\nu+1}+it$ parallel to the critical line $1/2+it$ for each interval\footnote{$\Delta_\lambda$ denotes the tolerance for the shift zeros of $Z_{\alpha^{(\lambda)}_\nu}(t)$ relative to zeros of $\cos\theta(t)$} $\Delta T_\lambda\pm\Delta_\lambda$
\begin{equation}
\label{function_sequence}
\begin{gathered}
Z_{\alpha^{(\lambda)}_\nu, 0}=Z_{\alpha^{(\lambda)}_{\nu+1}}\\
Z_{\alpha^{(\lambda)}_\nu, k}=Z_{\alpha^{(\lambda)}_\nu, k-1}-P_{\alpha^{(\lambda)}_\nu,k}
\end{gathered}
\end{equation}
\begin{equation}
\notag
\begin{split}
P_{\alpha^{(\lambda)}_\nu,k}=\frac{(-1)^k}{k!}\big\{\alpha^{(\lambda)}_{\nu+1}(\alpha^{(\lambda)}_{\nu+1}-1)(\alpha^{(\lambda)}_{\nu+1}-2)&...(\alpha^{(\lambda)}_{\nu+1}-k+1)\\
-\alpha^{(\lambda)}_{\nu}(\alpha^{(\lambda)}_{\nu}-1)(\alpha^{(\lambda)}_{\nu}-2)&...(\alpha^{(\lambda)}_{\nu}-k+1)\big\}Q_k
\end{split}
\end{equation}
\par 
Obviously, to prove theorem on the transformation of zeros $\cos\theta(t)$ into zeros of the Hardy function, we need to define the following sequence of functions
\begin{equation}
\notag
\begin{gathered}
\{f_k\}^\infty_{k=1}=\big\{\cos\theta\to Z_{\hat\alpha_\lambda}\to \{Z_{\alpha^{(\lambda)}_\nu, k}\}^\infty_{k=1}\to Z\big\}\\
\alpha_{\nu}=1/2,\ \alpha_{\nu+1}=\hat\alpha_\lambda
\end{gathered}
\end{equation}
\par 
To have possibility to decrease value $\max|P_{\alpha^{(\lambda)}_\nu,k}|$ between any neighboring zeros of $Z_{\alpha^{(\lambda)}_\nu, k-1}$ we fix the $\mu_\lambda$ values 
$$\alpha^{(\lambda)}_1<\alpha^{(\lambda)}_2<\alpha^{(\lambda)}_3<...<\alpha^{(\lambda)}_{\nu}<...<\alpha^{(\lambda)}_{\mu_\lambda}$$ 
between the values $\sigma=1/2$ and $\sigma=\hat\alpha_\lambda>A_\lambda$ for each interval $\Delta T_\lambda\pm\Delta_\lambda$ and define a sequence of functions

\begin{equation}
\notag
\begin{split}
\{f^{(\lambda)}_k\}^\infty_{k=1} &=\big\{\cos\theta\to Z_{\hat\alpha_\lambda}\to \{Z_{\alpha^{(\lambda)}_{\mu_\lambda}, k}\}^\infty_{k=1}\\
&\to Z_{\alpha^{(\lambda)}_{\mu_\lambda}}\to \{Z_{\alpha^{(\lambda)}_{\mu_\lambda-1}, k}\}^\infty_{k=1}\\
...&\to Z_{\alpha^{(\lambda)}_{\nu+1}}\to \{Z_{\alpha^{(\lambda)}_\nu, k}\}^\infty_{k=1}\to Z_{\alpha^{(\lambda)}_{\nu}}...\\
&\to Z_{\alpha^{(\lambda)}_3}\to  \{Z_{\alpha^{(\lambda)}_2, k}\}^\infty_{k=1}\\
&\to Z_{\alpha^{(\lambda)}_2}\to  \{Z^{(\lambda)}_{\alpha^{(\lambda)}_1, k}\}^\infty_{k=1}\to Z\big\}
\end{split}
\end{equation}
\par 
and build pairs of zeros $\{f^{(\lambda)}_k\}^\infty_{k=1}$ in the required sequence
\begin{equation}
\notag
\begin{split}
&\{\overline t^{(\lambda)}_n,\ \overline t^{(\lambda)}_{n+1}\}\to \{\delta_{\hat\alpha_\lambda,n},\ \delta_{\hat\alpha_\lambda,n+1}\}\to \{\delta_{\alpha^{(\lambda)}_{\mu_\lambda},k,n},\ \delta_{\alpha^{(\lambda)}_{\mu_\lambda},k,n+1}\}^\infty_{k=1}\\
&\to \{\delta_{\alpha^{(\lambda)}_{\mu_\lambda},n},\ \delta_{\alpha^{(\lambda)}_{\mu_\lambda},n+1}\}\to \{\delta_{\alpha^{(\lambda)}_{\mu_\lambda-1},k,n},\ \delta
_{\alpha^{(\lambda)}_{\mu_\lambda-1},k,n+1}\}^\infty_{k=1}\\
...&\to \{\delta_{\alpha^{(\lambda)}_{\nu+1},n},\ \delta_{\alpha^{(\lambda)}_{\nu+1},n+1}\}\to \{\delta_{\alpha^{(\lambda)}_\nu,k,n},\ \delta_{\alpha^{(\lambda)}_\nu,k,n+1}\}^\infty_{k=1}\\
&\to \{\delta_{\alpha^{(\lambda)}_\nu,n},\ \delta_{\alpha^{(\lambda)}_\nu,n+1}\}...\\
&\to \{\delta_{\alpha^{(\lambda)}_3,n},\ \delta_{\alpha^{(\lambda)}_3,n+1}\}\to \{\delta_{\alpha^{(\lambda)}_2,k,n},\ \delta_{\alpha^{(\lambda)}_2,k,n+1}\}^\infty_{k=1}\\
&\to \{\delta_{\alpha^{(\lambda)}_2,n},\ \delta_{\alpha^{(\lambda)}_2,n+1}\}\to \{\delta^{(\lambda)}_{\alpha_1,k,n},\ \delta^{(\lambda)}_{\alpha_1,k,n+1}\}^\infty_{k=1}\to \{\gamma^{(\lambda)}_n,\ \gamma^{(\lambda)}_{n+1}\}
\end{split}
\end{equation}
\par 
where
\begin{equation}
\label{condition_t_lambda}
\begin{gathered}
\cos\theta(\overline t^{(\lambda)}_n)=0,\ Z_{\hat\alpha_\lambda}(\delta_{\hat\alpha_\lambda,n})=0,\ Z_{\alpha^{(\lambda)}_{\nu}}(\delta_{\alpha^{(\lambda)}_{\nu},n})=0\\
Z_{\alpha^{(\lambda)}_\nu, k}(\delta_{\alpha^{(\lambda)}_\nu,k,n})=0,\ Z(\gamma^{(\lambda)}_n)=0
\end{gathered}
\end{equation}
\par 
Thus, we proved the first part of the theorem on the transformation of zeros $\cos\theta(t)$ into zeros of the Hardy function.
\par 
\section{Sequence of zeros \texorpdfstring{$Z_{\alpha^{(\lambda)}_\nu, k}$}{}}
Review the transformation of the Generalized Hardy function on any interval $\Delta T_\lambda\pm\Delta_\lambda$ for any two values of $\alpha^{(\lambda)}_\nu<\alpha^{(\lambda)}_{\nu+1}$
$$Z_{\alpha^{(\lambda)}_{\nu+1}}\to \{Z_{\alpha^{(\lambda)}_\nu, k}\}^\infty_{k=1}\to Z_{\alpha^{(\lambda)}_{\nu}}$$
and the corresponding transformation of zeros
$$\{\delta_{\alpha^{(\lambda)}_{\nu+1},n},\ \delta_{\alpha^{(\lambda)}_{\nu+1},n+1}\}\to \{\delta_{\alpha^{(\lambda)}_\nu,k,n},\ \delta_{\alpha^{(\lambda)}_\nu,k,n+1}\}^\infty_{k=1}\to \{\delta_{\alpha^{(\lambda)}_\nu,n},\ \delta_{\alpha^{(\lambda)}_\nu,n+1}\}$$
\par 
where
$$Z_{\alpha^{(\lambda)}_{\nu}}(\delta_{\alpha^{(\lambda)}_{\nu},n})=0,\ Z_{\alpha^{(\lambda)}_\nu, k}(\delta_{\alpha^{(\lambda)}_\nu,k,n})=0$$

\par 
Obviously the zeros are $\{\delta_{\alpha^{(\lambda)}_\nu,k,n},\ \delta_{\alpha^{(\lambda)}_\nu,k,n+1}\}$ can be found by solving a system of equations  (\ref{function_sequence}).
\par 
To prove the second part of the theorem, we reviewed different options of ranges between points where the function $Z_{\alpha^{(\lambda)}_\nu, k-1}$ changes sign and $P_{\alpha^{(\lambda)}_\nu,k}>0$ or $P_{\alpha^{(\lambda)}_\nu,k}<0$.
 there are four options of ranges $\{1,\ 2,\ 3,\ 4\}$ comparing the values of $P_ {\alpha^ {(\lambda)} _ \nu, k} $ for 
$$Z_{\alpha^{(\lambda)}_\nu, k-1}(t)>0,\ t\in (\delta_{\alpha^{(\lambda)}_\nu,k-1,n}, \delta_{\alpha^{(\lambda)}_\nu,k-1,n+1})$$
\par 
1) $P_{\alpha^{(\lambda)}_\nu,k}(\delta_{\alpha^{(\lambda)}_\nu,k-1,n})>0$ and $P_{\alpha^{(\lambda)}_\nu,k}(\delta_{\alpha^{(\lambda)}_\nu,k-1,n+1})>0$
\par 
2) $P_{\alpha^{(\lambda)}_\nu,k}(\delta_{\alpha^{(\lambda)}_\nu,k-1,n})<0$ and $P_{\alpha^{(\lambda)}_\nu,k}(\delta_{\alpha^{(\lambda)}_\nu,k-1,n+1})<0$
\par 
3) $P_{\alpha^{(\lambda)}_\nu,k}(\delta_{\alpha^{(\lambda)}_\nu,k-1,n})>0$ and $P_{\alpha^{(\lambda)}_\nu,k}(\delta_{\alpha^{(\lambda)}_\nu,k-1,n+1})<0$
\par 
4) $P_{\alpha^{(\lambda)}_\nu,k}(\delta_{\alpha^{(\lambda)}_\nu,k-1,n})<0$ and $P_{\alpha^{(\lambda)}_\nu,k}(\delta_{\alpha^{(\lambda)}_\nu,k-1,n+1})>0$
\par 
and the same four options of ranges $\{1',\ 2',\ 3',\ 4'\}$ for 
$$Z_{\alpha^{(\lambda)}_\nu, k-1}(t)<0,\ t\in (\delta_{\alpha^{(\lambda)}_\nu,k-1,n}, \delta_{\alpha^{(\lambda)}_\nu,k-1,n+1})$$
\par 
At points where $Z_{\alpha^{(\lambda)}_\nu, k-1}(\delta_{\alpha^{(\lambda)}_\nu,k-1,n})=0$ and $Z_{\alpha^{(\lambda)}_\nu, k-1}(\delta_{\alpha^{(\lambda)}_\nu,k-1,n+1})=0$ the function $Z_{\alpha^{(\lambda)}_\nu, k-1}$ changes the sign, thus the options of ranges $\{1,\ 2,\ 3,\ 4\}$ and $\{1',\ 2',\ 3',\ 4'\}$  can change each other in a restricted number of ways
\begin{equation}
\label{options}
\begin{gathered}
\xymatrix{
&&(1)&&(3)\\
(1)\ar[r]^A\ar[dr]&(1')\ar[r]^B\ar[ur]&(3)\ar[r]^B\ar[dr]&(4')\ar[r]^B\ar[ur]&(1)\\
(4)\ar[r]^A\ar[ur]&(3')\ar[r]^A\ar[dr]&(2)\ar[r]^B\ar[ur]&(2')\ar[r]^A\ar[dr]&(4)\\
&&(4)&&(2)
}
\end{gathered}
\end{equation}
\par 
Obviously, following (\ref{function_sequence}), we can choose $\Delta\alpha^{(\lambda)}_\nu=\alpha^{(\lambda)}_{\nu+1}-\alpha^{(\lambda)}_\nu$ so as to keep the condition $\max|Z_{\alpha^{(\lambda)}_{\nu,k-1}}|>\max|P_{\alpha^{(\lambda)}_\nu,k}|$ at each range $(\delta_{\alpha^{(\lambda)}_\nu,k-1,n}, \delta_{\alpha^{(\lambda)}_\nu,k-1,n+1})$.
\par 
Then the analysis of the sequence of the options of ranges each after other (\ref{options}) shows that they can be matched with two options of the sequence of zeros
\par 
A) $\delta_{\alpha^{(\lambda)}_\nu,k,n}<\delta_{\alpha^{(\lambda)}_\nu,k-1,n}$
\par 
B) $\delta_{\alpha^{(\lambda)}_\nu,k,n}>\delta_{\alpha^{(\lambda)}_\nu,k-1,n}$
\par 
Really
\begin{equation}
\notag
\begin{gathered}
\forall(\delta_{\alpha^{(\lambda)}_\nu,k-1,n}, \delta_{\alpha^{(\lambda)}_\nu,k-1,n+1})\\
\exists C_\lambda:\ \max|Z_{\alpha^{(\lambda)}_{\nu+1}}(t)|>C_\lambda,\
\exists D_\lambda>\max\big\{|Q_k(t)|\big\},\ t\in (\delta_{\alpha^{(\lambda)}_\nu,k-1,n}, \delta_{\alpha^{(\lambda)}_\nu,k-1,n+1})\\
\exists \Delta \alpha^{(\lambda)}_\nu=\alpha^{(\lambda)}_{\nu+1}-\alpha^{(\lambda)}_\nu:
\end{gathered}
\end{equation}
\begin{equation}
\notag
\begin{split}
C_\lambda>\big|\frac{(-1)^k}{k!}\big\{\alpha^{(\lambda)}_{\nu+1}(\alpha^{(\lambda)}_{\nu+1}-1)(\alpha^{(\lambda)}_{\nu+1}-2)&...(\alpha^{(\lambda)}_{\nu+1}-k+1)\\
-\alpha^{(\lambda)}_{\nu}(\alpha^{(\lambda)}_{\nu}-1)(\alpha^{(\lambda)}_{\nu}-2)&...(\alpha^{(\lambda)}_{\nu}-k+1)\big\}\big|D_\lambda
\end{split}
\end{equation}
\par 
Then for any range $(\delta_{\alpha^{(\lambda)}_\nu,k-1,n}, \delta_{\alpha^{(\lambda)}_\nu,k-1,n+1})$ we can define the sequence $$\alpha^{(\lambda)}_1<\alpha^{(\lambda)}_2<\alpha^{(\lambda)}_3<...<\alpha^{(\lambda)}_{\nu}<...<\alpha^{(\lambda)}_{\mu_\lambda}$$ 
between the values $\sigma=1/2$ and $\sigma=\hat\alpha_\lambda>A_\lambda$ based on the condition
\begin{equation}
\label{condition_c_lambda}
\begin{gathered}
\max|Z_{\alpha^{(\lambda)}_{\nu+1}}(t)|>C_\lambda>\max\big\{|P_{\alpha^{(\lambda)}_\nu,k}(t)|\big\}\\
t\in (\delta_{\alpha^{(\lambda)}_\nu,k-1,n}, \delta_{\alpha^{(\lambda)}_\nu,k-1,n+1})
\end{gathered}
\end{equation}
\par 
Take, for example, a sequence of ranges
\xymatrix{
(1)\ar[r]^A&(1')\ar[r]^B&(3)
}
\begin{equation}
\label{ranges_example}
\begin{split}
Z_{\alpha^{(\lambda)}_\nu, k-1}(t)>0,\ t\in (\delta_{\alpha^{(\lambda)}_\nu,k-1,n-1}, \delta_{\alpha^{(\lambda)}_\nu,k-1,n})&=(1)\\
Z_{\alpha^{(\lambda)}_\nu, k-1}(t)<0,\ t\in (\delta_{\alpha^{(\lambda)}_\nu,k-1,n}, \delta_{\alpha^{(\lambda)}_\nu,k-1,n+1})&=(1')\\
Z_{\alpha^{(\lambda)}_\nu, k-1}(t)>0,\ t\in (\delta_{\alpha^{(\lambda)}_\nu,k-1,n+1}, \delta_{\alpha^{(\lambda)}_\nu,k-1,n+2})&=(3)
\end{split}
\end{equation}
\par 
and write out all the comparisons at the points where the function $Z_{\alpha^{(\lambda)}_\nu, k-1}$ changes sign
\begin{equation}
\notag
\begin{split}
P_{\alpha^{(\lambda)}_\nu,k}(\delta_{\alpha^{(\lambda)}_\nu,k-1,n-1})&>0\\
P_{\alpha^{(\lambda)}_\nu,k}(\delta_{\alpha^{(\lambda)}_\nu,k-1,n})&>0\\
P_{\alpha^{(\lambda)}_\nu,k}(\delta_{\alpha^{(\lambda)}_\nu,k-1,n+1})&>0\\
P_{\alpha^{(\lambda)}_\nu,k}(\delta_{\alpha^{(\lambda)}_\nu,k-1,n+2})&<0
\end{split}
\end{equation}
\par 
Review the ranges of $(1)$ and $(3)$, where $Z_{\alpha^{(\lambda)}_\nu, k-1}$ and $P_{\alpha^{(\lambda)}_\nu,k}$ have the same sign, then it follows from the condition (\ref{condition_c_lambda}) that inside the ranges $(1)$ and $(3)$  $Z_{\alpha^{(\lambda)}_\nu, k-1}>P_{\alpha^{(\lambda)}_\nu,k}$, while at the boundary of the range $Z_{\alpha^{(\lambda)}_\nu, k-1}<P_{\alpha^{(\lambda)}_\nu,k}$, then the expression $Z_{\alpha^{(\lambda)}_\nu, k-1}-P_{\alpha^{(\lambda)}_\nu,k}$ in the range $(1)$ changes the sign from right to left, hence we get a comparison 
$$A:\ \delta_{\alpha^{(\lambda)}_\nu,k,n}<\delta_{\alpha^{(\lambda)}_\nu,k-1,n}$$
and in the range $(3)$ - on the contrary- from left to right, therefore, we get a comparison
$$B:\ \delta_{\alpha^{(\lambda)}_\nu,k,n+1}>\delta_{\alpha^{(\lambda)}_\nu,k-1,n+1}$$
similarly, other sequence of ranges options (\ref{options})can be reviewed.
\par 
We have reviewed why zeros $\delta_{\alpha^{(\lambda)}_\nu,k,n}$ and $\delta_{\alpha^{(\lambda)}_\nu,k,n+1}$ are real, now we review why they form monotonic sequence.
\par 
From the expression
\begin{equation}
\notag
\begin{split}
C_{\nu,k}=\frac{(-1)^k}{k!}\big\{\alpha^{(\lambda)}_{\nu+1}(\alpha^{(\lambda)}_{\nu+1}-1)(\alpha^{(\lambda)}_{\nu+1}-2)&...(\alpha^{(\lambda)}_{\nu+1}-k+1)\\
-\alpha^{(\lambda)}_{\nu}(\alpha^{(\lambda)}_{\nu}-1)(\alpha^{(\lambda)}_{\nu}-2)&...(\alpha^{(\lambda)}_{\nu}-k+1)\big\}
\end{split}
\end{equation}
\par 
it can be obtained trivial\footnote{dependency can be obtained by direct calculation for any $\Delta\alpha^{(\lambda)}_\nu=\alpha^{(\lambda)}_{\nu+1}-\alpha^{(\lambda)}_\nu$} that
\begin{equation}
\notag
\begin{gathered}
\forall\Delta \alpha^{(\lambda)}_\nu=\alpha^{(\lambda)}_{\nu+1}-\alpha^{(\lambda)}_\nu\\ 
\exists M_1: \forall k>M_1\\
|C_{\nu, k}|>|C_{\nu, k+1}|>|C_{\nu, k+2}|>...
\end{gathered}
\end{equation}
\par 
From the expression (\ref{generalized_hardy_function_4}) it is easy to see that the sequence $\{|Q_k|\}$ is bounded from above by the sequence of functions 
$$\Big\{m\Big(\frac{m-1}{m}\Big)^k\Big\}^\infty_{k=1}$$
which forms a limited monotonic decreasing sequence for a fixed $m$
\begin{equation}
\notag
m\Big(\frac{m-1}{m}\Big)>m\Big(\frac{m-1}{m}\Big)^2>m\Big(\frac{m-1}{m}\Big)^3>...>m\Big(\frac{m-1}{m}\Big)^k>...
\end{equation}
\par 
Then
\begin{equation}
\notag
\begin{gathered}
\exists M_2: \forall k>M_2\\
|Q_k|>|Q_{k+1}|>|Q_{k+2}|>...
\end{gathered}
\end{equation}
\par 
In addition, using the asymptotic dependence
\begin{equation}
\notag
\begin{split}
\Big(\frac{m-1}{m}\Big)^k&\sim{e^{-\frac{k}{m}}}\\
m&\to\infty\
\end{split}
\end{equation}
\par 
sequence limits exist and both are zero
\begin{equation}
\notag
\begin{gathered}
\lim_{k\to\infty}|Q_k|=\lim_{k\to\infty}m\Big(\frac{m-1}{m}\Big)^k=0\\
m= \Big\lfloor\frac{k_m|t|}{2\pi}\Big\rfloor,\ k_m=3,4,5...
\end{gathered}
\end{equation}
\par 
Then $|P_{\alpha^{(\lambda)}_\nu,k}|$ also form monotonically decreasing sequence
\begin{equation}
\notag
\begin{gathered}
\forall k>\max(M_1,\ M_2)\\
|P_{\alpha^{(\lambda)}_\nu,k}|>|P_{\alpha^{(\lambda)}_\nu,k+1}|>|P_{\alpha^{(\lambda)}_\nu,k+2}|>...\\
\lim_{k\to\infty}|P_{\alpha^{(\lambda)}_\nu,k}|=0
\end{gathered}
\end{equation}
\par 
Ranges analysis (\ref{ranges_example}) showed that the zeros of $Z_{\alpha^{(\lambda)}_\nu,k}$ are formed in the range where $Z_{\alpha^{(\lambda)}_\nu, k-1}$ and $P_{\alpha^{(\lambda)}_\nu,k}$ have the same sign.
\par 
Then from (\ref{function_sequence}) we obtain that  $|Z_{\alpha^{(\lambda)}_\nu,k}|$ form a monotonically decreasing sequence on this range
\begin{equation}
\label{decreasing_sequence}
\begin{gathered}
\forall k>\max(M_1,\ M_2)\\
|Z_{\alpha^{(\lambda)}_\nu,k}|>|Z_{\alpha^{(\lambda)}_\nu,k+1}|>|Z_{\alpha^{(\lambda)}_\nu,k+2}|>...\\
\end{gathered}
\end{equation}
\par 
Therefore, we can apply to the relation $|P_{\alpha^{(\lambda)}_\nu,k}|/|Z_{\alpha^{(\lambda)}_\nu,k-1}|$ the Stolz-Cesaro theorem on sequences \cite{CN}, a variant for $a_n\to 0, b_n\to 0, b_n-b_{n+1}>0$, i.e. $b_n$ is strictly decreasing
\begin{equation}
\notag
\lim_{n\to\infty}\frac{a_{n+1}-a_n}{b_{n+1}-b_n}=\lim_{n\to\infty}\frac{a_n}{b_n}=l
\end{equation}
\par 
Obviously we can write (\ref{function_sequence}) as
$$|P_{\alpha^{(\lambda)}_\nu,k}|=|Z_{\alpha^{(\lambda)}_\nu,k}|-|Z_{\alpha^{(\lambda)}_\nu,k-1}|$$
then
\begin{equation}
\notag
\begin{gathered}
\lim_{k\to\infty}\frac{|P_{\alpha^{(\lambda)}_\nu,k+1}|-|P_{\alpha^{(\lambda)}_\nu,k}|}{|Z_{\alpha^{(\lambda)}_\nu,k}|-|Z_{\alpha^{(\lambda)}_\nu,k-1}|}=\lim_{k\to\infty}\frac{|P_{\alpha^{(\lambda)}_\nu,k+1}|-|P_{\alpha^{(\lambda)}_\nu,k}|}{|P_{\alpha^{(\lambda)}_\nu,k}|}=\lim_{k\to\infty}\frac{|P_{\alpha^{(\lambda)}_\nu,k+1}|}{|P_{\alpha^{(\lambda)}_\nu,k}|}-\lim_{k\to\infty}\frac{|P_{\alpha^{(\lambda)}_\nu,k}|}{|P_{\alpha^{(\lambda)}_\nu,k}|}=1-1=0
\end{gathered}
\end{equation}
\par 
Therefore
\begin{equation}
\label{limit_4}
\lim_{k\to\infty}\frac{|P_{\alpha^{(\lambda)}_\nu,k}|}{|Z_{\alpha^{(\lambda)}_\nu,k-1}|}=\lim_{k\to\infty}\frac{|P_{\alpha^{(\lambda)}_\nu,k+1}|-|P_{\alpha^{(\lambda)}_\nu,k}|}{|Z_{\alpha^{(\lambda)}_\nu,k}|-|Z_{\alpha^{(\lambda)}_\nu,k-1}|}=0
\end{equation}
\par 
The result now easy follows from (\ref{decreasing_sequence}, \ref{limit_4}) in the case of $A:\ \delta_{\alpha^{(\lambda)}_\nu,k,n}<\delta_{\alpha^{(\lambda)}_\nu,k-1,n}$ zeros $\delta_{\alpha^{(\lambda)}_\nu,k,n}$ form a strictly decreasing sequence of real zeros
\begin{equation}
\notag
\begin{gathered}
\forall k>\max(M_1,\ M_2)\\
\delta_{\alpha^{(\lambda)}_\nu,k,n}>\delta_{\alpha^{(\lambda)}_\nu,k+1,n}>\delta_{\alpha^{(\lambda)}_\nu,k+2,n}>...\\
\end{gathered}
\end{equation}
\par 
Hence, zero $\delta_{\alpha^{(\lambda)}_\nu,n}$ is also real, as is the limit of this sequence
$$\lim_{k\to\infty}\delta_{\alpha^{(\lambda)}_\nu,k,n}=\delta_{\alpha^{(\lambda)}_\nu,n}$$
\par 
Similarly in the case of $B:\ \delta_{\alpha^{(\lambda)}_\nu,k,n}>\delta_{\alpha^{(\lambda)}_\nu,k-1,n}$ zeros $\delta_{\alpha^{(\lambda)}_\nu,k,n}$ form a strictly increasing sequence of real zeros
\begin{equation}
\notag
\begin{gathered}
\forall k>\max(M_1,\ M_2)\\
\delta_{\alpha^{(\lambda)}_\nu,k,n}<\delta_{\alpha^{(\lambda)}_\nu,k+1,n}<\delta_{\alpha^{(\lambda)}_\nu,k+2,n}<...\\
\end{gathered}
\end{equation}
\par 
Then, zero $\delta_{\alpha^{(\lambda)}_\nu,n}$ is also real, as is the limit of this sequence
$$\lim_{k\to\infty}\delta_{\alpha^{(\lambda)}_\nu,k,n}=\delta_{\alpha^{(\lambda)}_\nu,n}$$
\par 
Thus, by induction between values $\sigma=1/2$ and $\sigma=\hat\alpha_\lambda>A_\lambda$
$$\alpha^{(\lambda)}_1<\alpha^{(\lambda)}_2<\alpha^{(\lambda)}_3<...<\alpha^{(\lambda)}_{\nu}<...<\alpha^{(\lambda)}_{\mu_\lambda}$$ 
an distinct relationship has been established between the zeros of the function $\cos\theta(t)$ and  the zeros of the Hardy function on any interval $\Delta T_\lambda\pm\Delta_\lambda$.

\par 
Thus, we proved the second part of the theorem on the transformation of zeros $\cos\theta(t)$ into zeros of the Hardy function.
\par 
\section{Conclusion}
\par 
In the theory of the Zeta function \cite{VK,TI}, the Riemann-Siegel theta function appears twice:
\par 
the first time in connection with the Hardy function
$$\theta(t)=\frac{t}{2}\log\frac{t}{2\pi}-\frac{t}{2}-\frac{\pi}{8}+O(1/t)$$
\par 
the second time is in connection with the Riemann-Mangoldt formula
$$N(T)=\frac{1}{\pi}Im\ \int_C\frac{\xi'(s)}{\xi(s)}ds)=\frac{1}{\pi}\theta(T)+1+\frac{1}{\pi}Im\ \int_C\frac{\zeta'(s)}{\zeta(s)}ds)$$
\par 
where the contour $C$ is defined by segments $[2,\ 2+iT]$ and $[ 2+iT,\ 1/2+iT]$
\par 
Obviously, the number of zeros $\cos\theta(t)$ in the interval (0, T] can be calculated using the formula
$$N_\theta(T)=\frac{1}{\pi}\theta(T)$$
\par 
In 1918, Backlund  \cite{BL} refined the results on estimating the variation of the argument of the Zeta function and the residual term of the theta function
$$Q(T)=N(T)-\Big(\frac{T}{2\pi}\log\frac{T}{2\pi}-\frac{T}{2\pi}+\frac{7}{8}\Big)=P(T)+\frac{1}{\pi}R(T)$$
$$|Q(T)|< 0.137\log T+0.443\log\log T+4.350$$
\par 
Thus, Backlund established a tolerance for the deviation of the number of zeros of the Zeta function from the number of zeros of the function $\cos\theta(t)$, which is due, as we found out, to the offset of zeros $Z_{\alpha^{(\lambda)}_\nu}(t)$ relative to the zeros of $\cos\theta(t)$.
\par 
Before we proved that any pair of zeros of the Hardy function  $\gamma_n,\ \gamma_{n+1}$ corresponds to zeros of the function $\cos\theta(t)$ $\overline t_n,\ \overline t_{n+1}$ and zeros of Hardy functions $\gamma_n,\ \gamma_{n+1}$ are real, now we have established the correspondence between the number of zeros of the function $\cos\theta(t)$ and the Zeta function\footnote{as we remember, the Hardy function $Z(t)$ takes real values for real $t$ and its real zeros are zeros $\zeta(s)$ on the critical line $1/2 + it$}.
\par 
Thus, we proved the third part of the theorem on the transformation of zeros $\cos\theta(t)$ into zeros of the Hardy function.
\par 
So, there was reason not to accept the assumption made by Lehmer\cite{LD} in 1956 that the Hardy function can have a negative local maximum or a positive local minimum.

\bibliographystyle{unsrt}  


\end{document}